\documentclass[12pt,twoside]{amsart}

\usepackage{amssymb}

\usepackage[author-year]{amsrefs}

\usepackage{amsthm}

\theoremstyle{plain}
\newtheorem{lemma}{Lemma}
\newtheorem{theorem}{Theorem}

\begin{document}

\title[Functional limit theorems]{Functional limit theorems for
multiparameter fractional Brownian motion}

\author{Anatoliy Malyarenko}

\thanks{This work is supported in part by the Foundation for
Knowledge and Competence Development}

\date{\today}

\address{M\"{a}lardalen University, Box 883, SE-721 23
V\"{a}ster{\aa}s, Sweden}

\urladdr{http://www.mdh.se/ima/forskning/forskare/001-anatoliy.malyarenko.cvnotes.shtml}

\email{anatoliy.malyarenko@mdh.se}



\begin{abstract}
We prove a general functional limit theorem for multiparameter
fractional Brownian motion. The functional law of the iterated
logarithm, functional L\'{e}vy's modulus of continuity and many
other results are its particular cases. Applications to
approximation theory are discussed.
\end{abstract}

\maketitle

\section{Introduction}

Let $B(t)=B(t,\omega)$, $t\geq 0$, $\omega\in\Omega$ be the
Brownian motion on the probability space
$(\Omega,\mathfrak{F},\mathsf{P})$. The \emph{law of the iterated
logarithm} \cite{Khi} states that
\[
\mathsf{P}\left\{\omega\colon\limsup_{t\to\infty}\frac{B(t,\omega)}{\sqrt{2t\log\log
t}}=1\right\}=1.
\]
We abbreviate this as
\begin{equation}\label{lil}
\limsup_{t\to\infty}\frac{B(t)}{\sqrt{2t\log\log
t}}=1\quad\mathsf{P}-\text{a. s.},
\end{equation}
where a. s. stands for almost surely.

The functional counterpart to the law of the iterated logarithm
was discovered by \cite{Str}. Let $C[0,1]$ be the Banach space of
all continuous functions $f\colon [0,1]\mapsto\mathbb{R}$ with the
uniform topology generated by the maximum norm
\[
\|f\|_{\infty}=\max_{t\in[0,1]}|f(t)|.
\]
Let $\mathcal{H}_B$ be the Hilbert space of all absolutely
continuous functions $f\colon [0,1]\mapsto\mathbb{R}$ with
$f(0)=0$ and finite \emph{Strassen's norm}
\[
\|f\|_S=\left(\int^1_0(f'(t))^2\,dt\right)^{1/2}.
\]
The centred unit ball $\mathcal{K}_B$ of the space $\mathcal{H}_B$
\[
\mathcal{K}_B=\{\,f\in\mathcal{H}_B\colon\|f\|_S\leq 1\,\}
\]
is called \emph{Strassen's ball}. Define
\[
\mathcal{S}=\left\{\,\eta_u(t)=\frac{B(ut)}{\sqrt{2u\log\log
u}}\colon u>e\,\right\}\subset C[0,1].
\]
The \emph{functional law of the iterated logarithm} states that,
in the uniform topology, the set of $\mathsf{P}$-a. s. limit
points of $\mathcal{S}$ as $u\to\infty$ is Strassen's
ball~$\mathcal{K}_B$. It follows that for any continuous
functional $F\colon C[0,1]\mapsto\mathbb{R}$
\begin{equation}\label{cf}
\limsup_{u\to\infty}F(\eta_u(t))=\sup_{f\in\mathcal{K}_B}F(f)\quad\mathsf{P}-\text{a.
s.}
\end{equation}
In particular, for $F(f)=f(1)$ the supremum
$\sup_{f\in\mathcal{K}_B}f(1)$ is equal to $1$ and attained on the
function $f(t)=t$. Therefore \eqref{cf} transforms into
\eqref{lil}.

Another interesting ordinary limit theorem is \emph{L\'{e}vy's
modulus of continuity} \cite{Lev}. It states that
\begin{equation}\label{lmc}
\limsup_{u\downarrow
0}\sup_{t\in[0,1]}\frac{|B(t+u)-B(t)|}{\sqrt{2u\log
u^{-1}}}=1\quad\mathsf{P}-\text{a. s.}
\end{equation}

The corresponding functional counterpart was discovered by
\cite{Mue}. Define
\[
\mathcal{S}(u)=\left\{\,\eta_s(t)=\frac{B(s+ut)-B(s)}{\sqrt{2u\log
u^{-1}}}\colon 0\leq s\leq 1-u\,\right\}\subset C[0,1].
\]
Then, in the uniform topology, the set of $\mathsf{P}$-a. s. limit
points of $\mathcal{S}(u)$ as $u\downarrow 0$ is Strassen's ball
$\mathcal{K}_B$. L\'{e}vy's modulus of continuity \eqref{lmc}
follows from its functional counterpart in the same way as the law
of the iterated logarithm \eqref{lil} follows from Strassen's law.

In fact, \cite{Mue} contains a general functional limit theorem
that includes the functional law of the iterated logaritm, the
functional L\'{e}vy's modulus of continuity and many other results
as particular cases. Our aim is to prove the analogue of the
results of \cite{Mue} for the \emph{multiparameter fractional
Brownian motion}. This is the separable centred Gaussian random
field $\xi(\mathbf{x})$ on the space $\mathbb{R}^N$ with the
covariance function
\begin{equation}\label{cov}
\begin{aligned}
R(\mathbf{x},\mathbf{y})&=\mathsf{E}\xi(\mathbf{x})\xi(\mathbf{y})\\
&=\frac{1}{2}(\|\mathbf{x}\|^{2H}+\|\mathbf{y}\|^{2H}-\|\mathbf{x}-\mathbf{y}\|^{2H}),
\end{aligned}
\end{equation}
where $\|\boldsymbol{\cdot}\|$ denotes the usual Euclidean norm on
the space $\mathbb{R}^N$. The parameter $H\in(0,1)$ is called the
\emph{Hurst parameter}. In particular, for $N=1$ and $H=1/2$ the
multiparameter fractional Brownian motion becomes
\[
\xi(t)=
\begin{cases}
B_1(t),&t\geq 0,\\
B_2(t),&t<0,
\end{cases}
\]
where $B_1(t)$ and $B_2(t)$ are two independent copies of the
Brownian motion.

In Section~\ref{formulation} we formulate our results. They are
proved in Section~\ref{proofs}. Examples and applications are
discussed in Section~\ref{examples}.

\section{Formulation of results}\label{formulation}

In what follows, we write $\xi_1(x)\stackrel{d}{=}\xi_2(x)$, if
two random functions $\xi_1(x)$ and $\xi_2(x)$ are defined on the
same space $X$ and have the same finite-dimensional distributions.
We denote by $O(N)$ the group of all orthogonal matrices on the
space $\mathbb{R}^N$.

\begin{lemma}\label{l1}
The multiparameter fractional Brownian motion has the next
properties.
\begin{enumerate}
\item\label{pro1} It has homogeneous increments, i.e., for any
$\mathbf{y}\in\mathbb{R}^N$
\begin{equation}\label{prop1}
\xi(\mathbf{x}+\mathbf{y})-\xi(\mathbf{y})\stackrel{d}{=}\xi(\mathbf{x}),
\end{equation}
\item\label{pro2} It is self-similar, i.e., for any
$u\in\mathbb{R}$
\begin{equation}\label{prop2}
\xi(u\mathbf{x})\stackrel{d}{=}u^H\xi(\mathbf{x}).
\end{equation}
\item\label{pro3} It is isotropic, i.e., for any $g\in O(N)$
\begin{equation}\label{prop3}
\xi(g\mathbf{x})\stackrel{d}{=}\xi(\mathbf{x}).
\end{equation}
\end{enumerate}
\end{lemma}

Property~\ref{pro3} prompts us to use the $O(N)$-invariant closed
unit ball of the space~$\mathbb{R}^N$:
\[
\mathcal{B}=\{\,\mathbf{x}\in\mathbb{R}^N\colon\|\mathbf{x}\|\leq
1\,\}
\]
(but \emph{not} the cube $[0,1]^N$, which is not
$O(N)$-invariant!) in the case of the multiparameter fractional
Brownian motion instead of the interval $[0,1]$ in the case of the
ordinary Brownian motion.

First of all, we need to describe Hilbert space
$\mathcal{H}_{\xi}$ and its closed unit ball $\mathcal{K}_{\xi}$.
In the case of the Brownian motion, the space $\mathcal{H}_B$ can
be characterised as the reproducing kernel Hilbert space for the
Brownian motion, or as the set of all admissible shifts of the
Gaussian measure $\mu_B$ on the space $C[0,1]$ that corresponds to
the Brownian motion, or as the kernel of the measure $\mu_B$
\cite{Lif}. In order to describe $\mathcal{H}_{\xi}$ we need to
introduce some notations.

Let $r,\varphi,\vartheta_1,\dots,\vartheta_{N-2}$ be the spherical
coordinates in $\mathcal{B}$. The set of \emph{spherical
harmonics} $S^l_m(\varphi,\vartheta_1,\dots,\vartheta_{N-2})$
forms the orthonormal basis in the Hilbert space $L^2(S^{N-1},dS)$
of all square integrable functions on the unit sphere $S^{N-1}$
with respect to the Lebesgue measure
\[
dS=\sin\vartheta_1\sin^2\vartheta_2\dots\sin^{N-2}\vartheta_{N-2}\,d\varphi
\,d\vartheta_1\,d\vartheta_2\dots\,d\vartheta_{N-2}.
\]
Here $m\geq 0$ and $1\leq l\leq h(m,N)$, where
\[
h(m,N)=\frac{(2m+N-2)(m+N-3)!}{(N-2)!m!}.
\]

Let $\delta^k_j$ denotes the Kronecker's symbol. Let
${}_pF_q(a_1,\dots,a_p;b_1,\dots,b_q;z)$ denotes the
hypergeometric function. Let
\[
\lambda_{m1}\geq\lambda_{m2}\geq\dots\geq\lambda_{mn}\geq\dots>0
\]
be the sequence of all eigenvalues (with multiplicities) of the
positive definite kernel
\begin{equation}\label{bm}
\begin{aligned}
b_m(r,s)&=\frac{\pi^{N/2}}{\Gamma(N/2+m)}\left[(r^{2H}+s^{2H})\delta^0_m
-\frac{\Gamma(m-H)}{\Gamma(-H)}(rs)^m(r+s)^{2(H-m)}\right.\\
&\quad\times\left.{}_2F_1\left(m+(N-2)/2,m-H;2m+N-1;\frac{4rs}{(r+s)^2}\right)\right]
\end{aligned}
\end{equation}
in the Hilbert space $L^2([0,1],dr)$. Let $\psi_{mn}(r)$ be the
eigenbasis of the kernel $b_m(r,s)$. The set
\begin{equation}\label{basis}
\{\,\psi_{mn}(r)S^l_m(\varphi,\vartheta_1,\dots,\vartheta_{N-2})\colon
m\geq 0,n\geq 1,1\leq l\leq h(m,N)\,\}
\end{equation}
forms the orthonormal basis in the Hilbert space
$L^2(\mathcal{B},dr\,dS)$. For any $f\in C(\mathcal{B})$, let
$f^l_{mn}$ be the Fourier coefficients of $f$ with respect to the
basis~\eqref{basis}:
\[
f^l_{mn}=\int_{S^{N-1}}\int^1_0f(r,\varphi,\vartheta_1,\dots,\vartheta_{N-2})
\psi_{mn}(r)S^l_m(\varphi,\vartheta_1,\dots,\vartheta_{N-2})\,dr\,dS.
\]

\begin{lemma}\label{l2}
The reproducing kernel Hilbert space $\mathcal{H}_{\xi}$ of the
multiparameter fractional Brownian motion $\xi(\mathbf{x})$,
$\mathbf{x}\in\mathcal{B}$ consists of all functions $f\in
C(\mathcal{B})$ with $f(\mathbf{0})=0$ that satisfy the condition
\[
\|f\|^2_S=\sum^{\infty}_{m=0}\sum^{\infty}_{n=1}\sum^{h(m,N)}_{l=1}
\frac{(f^l_{mn})^2}{\lambda_{mn}}<\infty.
\]
The scalar product in the space $\mathcal{H}_{\xi}$ is defined as
\[
(f,g)_S=\sum^{\infty}_{m=0}\sum^{\infty}_{n=1}\sum^{h(m,N)}_{l=1}
\frac{f^l_{mn}g^l_{mn}}{\lambda_{mn}}.
\]
\end{lemma}

Therefore Strassen's ball is described as
\[
\mathcal{K}_{\xi}=\{\,f\in\mathcal{H}_{\xi}\colon\|f\|^2_S\leq
1\,\}.
\]
In what follows we write $\mathcal{K}$ instead of
$\mathcal{K}_{\xi}$.

Let $t_0$ be a real number. Let for every $t\geq t_0$ there exists
a non-empty set of indices $\mathcal{J}(t)$. Let every element
$j\in\mathcal{J}(t)$ defines the vector
$\mathbf{y}_j\in\mathbb{R}^N$ and the positive real number $u_j$.
Let $R_r(\mathbf{y}_j,u_j)$ be the cylinder
\[
R_r(\mathbf{y}_j,u_j)=\{\,(\mathbf{y},u)\colon
\|\mathbf{y}-\mathbf{y}_j\|\leq ru_j,e^{-r}u_j\leq u\leq
e^ru_j\,\},\qquad r>0.
\]

Now we define the function $F_r(t)$. In words: this is the volume
of the union of all cylinders $R_r(\mathbf{y}_j,u_j)$ that are
defined before the moment $t$, with respect to the measure
$u^{-N-1}\,d\mathbf{y}\,du$. Formally,
\[
F_r(t)=\int_{\cup_{t_0\leq v\leq t}\cup_{j\in\mathcal{J}(v)}
R_r(\mathbf{y}_j,u_j)}u^{-N-1}\,d\mathbf{y}\,du.
\]

Finally, we define
\[
\mathcal{P}(t)=\{\,(\mathbf{y}_j,u_j)\colon j\in\mathcal{J}(t)\,\}
\]
and the \emph{cloud of normed increments}
\[
\mathcal{S}(t)=\left\{\,\eta(\mathbf{x})=
\frac{\xi(\mathbf{y}+u\mathbf{x})-\xi(\mathbf{y})}{\sqrt{2h(t)}u^H}\colon
(\mathbf{y},u)\in\mathcal{P}(t)\,\right\}\subset C(\mathcal{B}).
\]

\begin{theorem}\label{th1}
Let the function $h(t)\colon[t_0,\infty)\mapsto\mathbb{R}$
satisfies the next conditions:
\begin{enumerate}
\item\label{cond1} $h(t)$ is increasing and
$\displaystyle\lim_{t\to\infty}h(t)=\infty$.
\item\label{cond2} The integral
$\displaystyle\int^{\infty}_{t_0}e^{-ah(t)}\,dF_1(t)$ converges
for $a>1$ and diverges for $a<1$.
\end{enumerate}
Then, in the uniform topology, the set of $\mathsf{P}$-a. s. limit
points of the cloud of increments $\mathcal{S}(t)$ as $t\to\infty$
is Strassen's ball~$\mathcal{K}$.
\end{theorem}

\section{Proofs}\label{proofs}

\subsection{Proof of Lemmas \ref{l1} and \ref{l2}}

\begin{proof}[Proof of Lemma~\ref{l1}]
It is enough to calculate the covariance functions of both hand
sides in equations~\eqref{prop1}--\eqref{prop3}. Calculations are
straightforward.
\end{proof}

\begin{proof}[Proof of Lemma~\ref{l2}]
The covariance function \eqref{cov} can be written as function of
three variables:
\begin{equation}\label{var3}
R(r,s,t)=\frac{1}{2}(r^{2H}+s^{2H}-(r^2+s^2-2rst)^H),\qquad r\geq
0,s\geq 0,-1\leq t\leq 1,
\end{equation}
where $r=\|\mathbf{x}\|$, $s=\|\mathbf{y}\|$, and $t$ is the
cosine of the angle between vectors $\mathbf{x}$ and $\mathbf{y}$.
According to the general theory of isotropic random fields
\cite{Yad}, the multiparameter fractional Brownian motion can be
written as
\begin{equation}\label{isoseries}
\xi(r,\varphi,\vartheta_1,\dots,\vartheta_{N-2})=
\sum^{\infty}_{m=0}\sum^{h(m,N)}_{l=1}\xi^l_m(r)S^l_m(\varphi,\vartheta_1,\dots,\vartheta_{N-2}),
\end{equation}
where $\xi^l_m(r)$ is the sequence of independent centred Gaussian
processes on $[0,\infty)$ with the covariance functions
\begin{equation}\label{covproc}
\mathsf{E}\xi^l_m(r)\xi^l_m(s)=\frac{2\pi^{(N-1)/2}}
{\Gamma((N-1)/2)C^{(N-2)/2}_m(1)}\int^1_{-1}R(r,s,t)C^{(N-2)/2}_m(t)
(1-t^2)^{(N-3)/2}\,dt,
\end{equation}
and $C^{(N-2)/2}_m(t)$ are Gegenbauer's polynomials. Denote by
$b_m(r,s)$ the covariance function \eqref{covproc}. We will prove
that $b_m(r,s)$ is expressed as \eqref{bm}.

It follows from \eqref{var3} and \eqref{covproc} that the
covariance function $b_m(r,s)$ can be written as the difference of
two integrals:
\begin{equation}\label{differ}
\begin{aligned}
b_m(r,s)&=\frac{\pi^{(N-1)/2}}{\Gamma((N-1)/2)C^{(N-2)/2}_m(1)}
(r^{2H}+s^{2H})\int^1_{-1}C^{(N-2)/2}_m(t)(1-t^2)^{(N-3)/2}\,dt\\
&\quad-\frac{\pi^{(N-1)/2}}{\Gamma((N-1)/2)C^{(N-2)/2}_m(1)}
\int^1_{-1}(r^2+s^2-2rst)^HC^{(N-2)/2}_m(t) (1-t^2)^{(N-3)/2}\,dt.
\end{aligned}
\end{equation}
The first integral is non-zero if and only if $m=0$ \cite{Vil}. It
follows that the first term in \eqref{differ} is equal to
\begin{align*}
&\delta^m_0\frac{\pi^{(N-1)/2}}{\Gamma((N-1)/2)}(r^{2H}+s^{2H})
\int^1_{-1}\frac{C^{(N-2)/2}_0(t)}{C^{(N-2)/2}_0(1)}(1-t^2)^{(N-3)/2}\,dt\\
&\quad=\frac{\delta^m_0\pi^{(N-1)/2}}{\Gamma((N-1)/2)}(r^{2H}+s^{2H})
\int^1_{-1}(1-t^2)^{(N-3)/2}\,dt\\
&\quad=\frac{\delta^m_0\pi^{N/2}}{\Gamma(N/2)}(r^{2H}+s^{2H}).
\end{align*}
Here we used formula 2.2.3.1 from \cite{Pru}.

Rewrite the second term as
\begin{align*}
&-\frac{\pi^{(N-1)/2}m!(N-3)!}{\Gamma((N-1)/2)(m+N-3)!}(2rs)^H
\lim_{\alpha\to(N-1)/2}\int^1_{-1}\left(\frac{r^2+s^2}{2rs}-t\right)^H\times\\
&\quad C^{(N-2)/2}_m(t)(1+t)^{\alpha-1}(1-t)^{(N-3)/2}\,dt.
\end{align*}
Using formula 2.21.4.15 from \cite{Pru2}, we can express this
limit as
\begin{align*}
&-\frac{(-1)^m2^{N-2}\pi^{(N-1)/2}\Gamma((N-1)/2)(m-1)!(r+s)^{2H}}{(m+N-2)!}
\lim_{\alpha\to(N-1)/2}\frac{\Gamma(\alpha-(N-3)/2-m)}{\Gamma((N-1)/2-\alpha)}\\
&\quad\times\lim_{\alpha\to(N-1)/2}\frac{{}_3F_2\left(\frac{N-1}{2},-H,1;
\alpha-\frac{N-3}{2}-m,m+N-1;\frac{4rs}{(r+s)^2}\right)}{\Gamma(\alpha-\frac{N-3}{2}-m)}.
\end{align*}
The first limit is calculated as
\begin{align*}
\lim_{\alpha\to(N-1)/2}\frac{\Gamma(\alpha-(N-3)/2-m)}{\Gamma((N-1)/2-\alpha)}
&=\lim_{\beta\to 0}\frac{\Gamma(-\beta-m+1)}{\Gamma(\beta)}\\
&=\lim_{\beta\to
0}\frac{(-1)^{m-1}\Gamma(-\beta)}{(1+\beta)(2+\beta)\dots(m-1+\beta)\Gamma(\beta)}\\
&=\frac{(-1)^{m-1}}{(m-1)!}\lim_{\beta\to
0}\frac{\Gamma(-\beta)}{\Gamma(\beta)}\\
&=\frac{(-1)^{m-1}}{(m-1)!}.
\end{align*}
For the second limit we use formula 7.2.3.6 from \cite{Pru3}:
\begin{align*}
&\lim_{\alpha\to(N-1)/2}\frac{{}_3F_2\left(\frac{N-1}{2},-H,1;
\alpha-\frac{N-3}{2}-m,m+N-1;\frac{4rs}{(r+s)^2}\right)}{\Gamma(\alpha-\frac{N-3}{2}-m)}\\
&\quad=\frac{(4rs)^m\Gamma((N-1)/2+m)\Gamma(m-H)(m+N-2)!}
{(r+s)^{2m}\Gamma((N-1)/2)\Gamma(-H)(2m+N-2)!}\\
&\qquad\times{}_2F_1\left(m+(N-2)/2,m-H;2m+N-1;\frac{4rs}{(r+s)^2}\right).
\end{align*}
Collecting all terms together, we obtain \eqref{bm}.

By Mercer's theorem, function $b_m(r,s)$ may be written as
uniformly and absolutely convergent series
\[
b_m(r,s)=\sum^{\infty}_{n=1}\lambda_{mn}\psi_n(r)\psi_n(s),\qquad
r,s\in[0,1].
\]
It follows that the random process $\xi^l_m(r)$ has the form
\[
\xi^l_m(r)=\sum^{\infty}_{n=1}\sqrt{\lambda_{mn}}\xi^l_{mn}\psi_n(r),\qquad
r\in[0,1],
\]
where $\xi^l_{mn}$ are independent standard normal random
variables. Substituting this representation to \eqref{isoseries},
we obtain:
\begin{equation}\label{local}
\xi(r,\varphi,\vartheta_1,\dots,\vartheta_{N-2})=
\sum^{\infty}_{m=0}\sum^{h(m,N)}_{l=1}\sum^{\infty}_{n=1}\sqrt{\lambda_{mn}}
\xi^l_{mn}\psi_n(r)S^l_m(\varphi,\vartheta_1,\dots,\vartheta_{N-2}).
\end{equation}
We call \eqref{local} the \emph{local spectral representation} of
the multiparameter fractional Brownian motion, because it is valid
only for $r\in[0,1]$, i.e., in $\mathcal{B}$.

Now Lemma~\ref{l2} follows from \eqref{local} and the general
theory of Gaussian measures \cite{Lif}.
\end{proof}

\subsection{Asymptotic relations that are equivalent to Theorem~\ref{th1}}

In this subsection, we formulate two asymptotic relations and
prove that they are equivalent to Theorem~\ref{th1}.

\begin{lemma}\label{l3}
The cloud of increments $\mathcal{S}(t)$ is $\mathsf{P}$-a. s.
almost inside $\mathcal{K}$ , i.e.,
\begin{equation}\label{attract}
\lim_{t\to\infty}\sup_{\eta\in\mathcal{S}(t)}\inf_{f\in\mathcal{K}}
\|\eta-f\|_{\infty}=0\quad\mathsf{P}-\text{a. s.}
\end{equation}
\end{lemma}

\begin{lemma}\label{l4}
Any neighbourhood of any element $f\in\mathcal{K}$ is caught by
the cloud of increments $\mathcal{S}(t)$ infinitely often, i.e.,
\begin{equation}\label{enter}
\sup_{f\in\mathcal{K}}\liminf_{t\to\infty}\inf_{\eta\in\mathcal{S}(t)}\|\eta-f\|_{\infty}=0.
\end{equation}
\end{lemma}

It is obvious that \eqref{attract} and \eqref{enter} follow from
Theorem~\ref{th1}.

Conversely, on the one hand, it follows from \eqref{attract} that
the set of $\mathsf{P}$-a.s. limit points of $\mathcal{S}(t)$
contains in the closure of $\mathcal{K}$. On the other hand, it
follows from \eqref{enter} that $\mathcal{K}$ contains in the set
of $\mathsf{P}$-a.s. limit points of $\mathcal{S}(t)$. According
to general theory \cite{Lif}, $\mathcal{K}$ is compact. Therefore
it is closed, and we are done.

\subsection{Construction of the auxiliary sequences}

We divide the set $\mathbb{R}^N\times(0,\infty)$ onto
parallelepipeds
\[
R_{\mathbf{k}p}=\{\,(\mathbf{y},u)\colon k_jre^{pr}\leq y_j\leq
(k_j+1)re^{pr}\quad\text{for}\quad 1\leq j\leq N,e^{pr}\leq u\leq
e^{(p+1)r}\,\},
\]
where $\mathbf{k}\in\mathbb{Z}^N$ and $p\in\mathbb{Z}$. The next
Lemma describes the most important property of the parallelepipeds
$R_{\mathbf{k}p}$.

\begin{lemma}\label{l5}
For any $t\in[t_0,\infty)$ the union of all cylinders
$R_r(\mathbf{y}_j,u_j)$ that are defined before the moment $t$
contains in the union of finitely many parallelepipeds
$R_{\mathbf{k}p}$.
\end{lemma}

\begin{proof}
It follows from condition~\ref{cond2} of Theorem~\ref{th1} that
for any $t\in[t_0,\infty)$ the volume of all cylinders
$R_r(\mathbf{y}_j,u_j)$ that are defined before the moment $t$
with respect to the measure $u^{-N-1}\,d\mathbf{y}\,du$ is finite.
So it is enough to prove that the volume of any parallelepiped
$R_{\mathbf{k}p}$ with respect to the above mentioned measure is
also finite. We have
\begin{align*}
\int_{R_{\mathbf{k}p}}u^{-N-1}\,d\mathbf{y}\,du&\sim
\frac{r^Ne^{Npr}[e^{(p+1)r}-e^{pr}]}{e^{(N+1)pr}}\\
&\sim r^N(e^r-1)\\
&\sim r^{N+1}\qquad (r\downarrow 0).
\end{align*}
Here and in what follows we write $f(r)\sim g(r)$ ($r\downarrow
0$) if
\[
\lim_{r\downarrow 0}\frac{f(r)}{g(r)}=1.
\]
\end{proof}

\begin{lemma}\label{l6}
There exist the sequence of real numbers $t_q$ and the sequence of
parallelepipeds $R_{\mathbf{k}_qp_q}$, $q\geq 0$, that satisfy the
next conditions.
\begin{enumerate}
\item\label{l6item1} For any $q\geq 0$ and for any $\varepsilon>0$
there exists a real number $t\in(t_q,t_q+\varepsilon)$ such that
\[
\mathcal{P}(t)\cap R_{\mathbf{k}_qp_q}\neq\varnothing.
\]
\item\label{l6item2} If $r<2/\sqrt{N}$ and $a>1$, then
\[
\sum^{\infty}_{q=0}\exp(-ah(t_q))<\infty.
\]
\end{enumerate}
\end{lemma}

\begin{proof}
We use mathematical induction.

The real number $t_0$ is already constructed (it is involved in
the formulation of Theorem~1). According to Lemma~\ref{l5}, the
union of all cylinders $R_r(\mathbf{y}_j,u_j)$ that are defined
before the moment $t_0+1$, contains in the union of finitely many
parallelepipeds $R_{\mathbf{k}p}$. Therefore there exists a
parallelepiped $R_{\mathbf{k}_0p_0}$ which intersects with
infinitely many sets from the sequence $\mathcal{P}(t_0+1)$,
$\mathcal{P}(t_0+1/2)$, \dots, $\mathcal{P}(t_0+1/n)$, \dots.

Assume that the real numbers $t_0$, $t_1$, \dots, $t_q$, and the
parallelepipeds $R_{\mathbf{k}_0p_0}$, $R_{\mathbf{k}_1p_1}$,
\dots, $R_{\mathbf{k}_qp_q}$ are already constructed. Define
\[
t_{q+1}=\inf\{\,t>t_q\colon\mathcal{P}(t)\nsubseteq
R_{\mathbf{k}_0p_0}\cup R_{\mathbf{k}_1p_1}\cup\dots\cup
R_{\mathbf{k}_qp_q}\,\}.
\]
Parallelepiped $R_{\mathbf{k}_{q+1}p_{q+1}}$ is defined as a
parallelepiped that intersects with infinitely many sets from the
sequence $\mathcal{P}(t_q+1)$, $\mathcal{P}(t_q+1/2)$, \dots,
$\mathcal{P}(t_q+1/n)$, \dots. It means that
condition~\ref{l6item1} is satisfied.

In order to prove condition~\ref{l6item2}, define the function
$F'_r(t)$ as the volume of the union of all the parallelepipeds
$R_{\mathbf{k}p}$, that intersect with at least one set
$\mathcal{P}(v)$ for $v\in[t_0,t]$, with respect to the measure
$u^{-N-1}\,d\mathbf{y}\,du$. Formally,
\[
F'_r(t)=\int_{\cup_{(\mathbf{k},p)\in\mathbb{Z}^{N+1}\colon
R_{\mathbf{k}p}\cap(\cup_{t_0\leq v\leq
t}\mathcal{P}(v))\neq\varnothing}R_{\mathbf{k}p}}u^{-N-1}\,d\mathbf{y}\,du.
\]
The length of a side of a cube, which is inscribed in the ball of
radius $1$ in the space $\mathbb{R}^N$, is equal to $2/\sqrt{N}$.
It follows that if $r<2/\sqrt{N}$ and $(\mathbf{y},u)\in
R_{\mathbf{k}p}$, then $R_{\mathbf{k}p}\subset R_1(\mathbf{y},u)$.
Therefore we have $F'_r(t)\leq F_1(t)$ and
\begin{align*}
\sum^{\infty}_{q=0}\exp(-ah(t_q))&\sim\frac{1}{r^{N+1}}\int^{\infty}_{t_0}
\exp(-ah(t))\,dF'_r(t)\\
&\leq\frac{1}{r^{N+1}}\int^{\infty}_{t_0}\exp(-ah(t))\,dF_1(t)\\
&<\infty.
\end{align*}
\end{proof}

\subsection{Proof of Lemma~\ref{l3}}

Denote
\[
\eta_{\mathbf{y},u}(\mathbf{x})=\frac{\xi(\mathbf{y}+u\mathbf{x})-\xi(\mathbf{y})}
{\sqrt{2}u^H},\qquad\mathbf{x}\in\mathcal{B}.
\]
Using properties \ref{pro1} and \ref{pro2} (Lemma~\ref{l1}), we
obtain
\[
\eta_{\mathbf{y},u}(\mathbf{x})\stackrel{d}{=}\frac{1}{\sqrt{2}}\xi(\mathbf{x}).
\]

Let $(\mathbf{y}_q,u_q)$ be the centre of the parallelepiped
$R_{\mathbf{k}_qp_q}$. Let
$(r,\varphi,\vartheta_1,\dots,\vartheta_{N-2})$ be the spherical
coordinates of a point $\mathbf{x}\in\mathcal{B}$. Let
$(r_2,\varphi_2,\vartheta_{1,2},\dots,\vartheta_{N-2,2})$ be the
spherical coordinates of a point $\mathbf{z}\in\mathcal{B}$. Let
$b^{l'l''}_{m'n'm''n''qs}$ be the Fourier coefficients of the
function
\[
b_{qs}(\mathbf{x},\mathbf{z})=\mathsf{E}\eta_{\mathbf{y}_q,u_q}(\mathbf{x})
\eta_{\mathbf{y}_s,u_s}(\mathbf{z})
\]
with respect to the orthonormal basis
\[
\psi_{m'n'}(r)S^{l'}_{m'}(\varphi,\vartheta_1,\dots,\vartheta_{N-2})
\psi_{m''n''}(r_2)S^{l''}_{m''}(\varphi_2,\vartheta_{1,2},\dots,\vartheta_{N-2,2}).
\]
Let $\{\xi^{lq}_{mn}\}$, $q\geq 0$ be the sequence of series of
standard normal random variables, that are independent in every
series, with the next correlation between series:
\[
\mathsf{E}\xi^{l'q}_{m'n'}\xi^{l''s}_{m''n''}=\lambda^{-1/2}_{m'n'}
\lambda^{-1/2}_{m''n''}b^{l'l''}_{m'n'm''n''},\qquad q\neq s.
\]
Then we have:
\[
\eta_{\mathbf{y}_q,u_q}(\mathbf{x})=\frac{1}{\sqrt{2}}\sum^{\infty}_{m=0}
\sum^{h(m,N)}_{l=1}\sum^{\infty}_{n=1}\sqrt{\lambda_{mn}}
\xi^{lq}_{mn}\psi_n(r)S^l_m(\varphi,\vartheta_1,\dots,\vartheta_{N-2}).
\]

Let $m_0$ and $n_0$ be two natural numbers. Denote
\[
\eta^{(m_0,n_0)}_{\mathbf{y}_q,u_q}(\mathbf{x})=\frac{1}{\sqrt{2}}\sum^{m_0}_{m=0}
\sum^{h(m,N)}_{l=1}\sum^{n_0}_{n=1}\sqrt{\lambda_{mn}}
\xi^{lq}_{mn}\psi_n(r)S^l_m(\varphi,\vartheta_1,\dots,\vartheta_{N-2}).
\]
For any $\varepsilon>0$ consider the next three events:
\begin{align*}
A_{1q}(\varepsilon)&=\left\{\frac{\eta^{(m_0,n_0)}_{\mathbf{y}_q,u_q}}
{\sqrt{h(t_q)}}\notin\mathcal{K}_{\varepsilon/3}\right\},\\
A_{2q}(\varepsilon)&=\left\{\frac{\|\eta_{\mathbf{y}_q,u_q}-
\eta^{(m_0,n_0)}_{\mathbf{y}_q,u_q}\|_{\infty}}{\sqrt{h(t_q)}}>\frac{\varepsilon}{3}\right\},\\
A_{3q}(\varepsilon)&=\left\{\sup_{(\mathbf{y},u)\in
R_{\mathbf{k}_qp_q}}\frac{\|\eta_{\mathbf{y}_q,u_q}-
\eta_{\mathbf{y},u}\|_{\infty}}{\sqrt{h(t_q)}}>\frac{\varepsilon}{3}\right\},
\end{align*}
where $\mathcal{K}_{\varepsilon/3}$ denotes the
$\varepsilon/3$-neighbourhood of Strassen's ball $\mathcal{K}$ in
the space $C(\mathcal{B})$. To prove Lemma~\ref{l3}, it is enough
to prove, that for any $\varepsilon>0$ there exist natural numbers
$m_0=m_0(\varepsilon)$ and $n_0=n_0(\varepsilon)$ such that the
events $A_{1q}(\varepsilon)$, $A_{2q}(\varepsilon)$, and
$A_{3q}(\varepsilon)$ occur only finitely many times
$\mathsf{P}$-a. s. In other words,
\[
\mathsf{P}\left\{\limsup_{q\to\infty}A_{1q}(\varepsilon)\right\}=
\mathsf{P}\left\{\limsup_{q\to\infty}A_{2q}(\varepsilon)\right\}=
\mathsf{P}\left\{\limsup_{q\to\infty}A_{3q}(\varepsilon)\right\}=0.
\]
By Borel--Cantelli lemma, it is enough to prove that
\begin{subequations}\label{inequal}
\begin{align}
\sum^{\infty}_{q=1}\mathsf{P}\{A_{1q}(\varepsilon)\}&<\infty,\\
\sum^{\infty}_{q=1}\mathsf{P}\{A_{2q}(\varepsilon)\}&<\infty,\\
\sum^{\infty}_{q=1}\mathsf{P}\{A_{3q}(\varepsilon)\}&<\infty.
\end{align}
\end{subequations}

We prove (\ref{inequal}b) first. Denote
\[
\sigma^2_{m_0n_0}(\mathbf{x})=\mathsf{E}\left[\eta_{\mathbf{y}_q,u_q}
(\mathbf{x})-\eta^{(m_0,n_0)}_{\mathbf{y}_q,u_q}(\mathbf{x})\right]^2,
\qquad\sigma^2_{m_0n_0}=\max_{\mathbf{x}\in\mathcal{B}}
\sigma^2_{m_0n_0}(\mathbf{x}).
\]
Using the large deviations estimate \cite{Lif}*{Section~12, (11)},
we obtain
\[
\mathsf{P}\left\{\|\eta_{\mathbf{y}_q,u_q}-
\eta^{(m_0,n_0)}_{\mathbf{y}_q,u_q}\|_{\infty}>\frac{\varepsilon\sqrt{h(t_q)}}{3}
\right\}\leq\exp\left[-\frac{\varepsilon^2h(t_q)}{18\sigma^2_{m_0n_0}}+o
\left(\frac{\varepsilon\sqrt{h(t_q)}}{3}\right)\right].
\]
By Lemma~\ref{l6}, condition~\ref{l6item2}, it is sufficient to
prove that for any $\varepsilon>0$ there exist natural numbers
$m_0=m_0(\varepsilon)$ and $n_0=n_0(\varepsilon)$ such that, say,
\[
\frac{\varepsilon^2}{9\sigma^2_{m_0n_0}}<1.
\]
Denote
\[
\eta^{(m_0)}_{\mathbf{y}_q,u_q}(\mathbf{x})=\frac{1}{\sqrt{2}}\sum^{m_0}_{m=0}
\sum^{h(m,N)}_{l=1}\sum^{\infty}_{n=1}\sqrt{\lambda_{mn}}
\xi^{lq}_{mn}\psi_n(r)S^l_m(\varphi,\vartheta_1,\dots,\vartheta_{N-2})
\]
and
\[
\sigma^2_{m_0}(\mathbf{x})=\mathsf{E}\left[\eta_{\mathbf{y}_q,u_q}
(\mathbf{x})-\eta^{(m_0)}_{\mathbf{y}_q,u_q}(\mathbf{x})\right]^2,
\qquad\sigma^2_{m_0}=\max_{\mathbf{x}\in\mathcal{B}}
\sigma^2_{m_0}(\mathbf{x}).
\]
The sequence $\sigma^2_{m_0}(\mathbf{x})$ converges to zero as
$m_0\to\infty$ for all $\mathbf{x}\in\mathcal{B}$. Moreover,
$\sigma^2_{m_0}(\mathbf{x})\geq\sigma^2_{m_0+1}(\mathbf{x})$ for
all $\mathbf{x}\in\mathcal{B}$ and all natural $m_0$. Functions
$\sigma^2_{m_0}(\mathbf{x})$ are non-negative and continuous. By
Dini's theorem, the sequence $\sigma^2_{m_0}(\mathbf{x})$
converges to zero uniformly on $\mathcal{B}$, i.e.,
\[
\lim_{m_0\to\infty}\sigma^2_{m_0}=0,
\]
and we choose such an $m_0$, that for any $m>m_0$,
$\sigma^2_m<\varepsilon^2/18$.

In the same way, we can apply Dini's theorem to the sequence of
functions
\[
\mathsf{E}\left[\eta^{(m_0)}_{\mathbf{y}_q,u_q}(\mathbf{x})-
\eta^{(m_0,n_0)}_{\mathbf{y}_q,u_q}(\mathbf{x})\right]^2,\quad
n_0\geq 1
\]
and find such $n_0$ that
\[
\sup_{\mathbf{x}\in\mathcal{B}}\mathsf{E}\left[\eta^{(m_0)}_{\mathbf{y}_q,u_q}(\mathbf{x})-
\eta^{(m_0,n)}_{\mathbf{y}_q,u_q}(\mathbf{x})\right]^2\leq\varepsilon^2/18
\]
for all $n>n_0$. (\ref{inequal}b) is proved.

Now we prove (\ref{inequal}a). In what follows we denote by $C$ a
constant depending only on $N$, $H$ and that may vary at each
occurrence. Specific constants will be denote by $C_1$, $C_2$,
\dots.

Consider the finite-dimensional subspace $E$ of the space
$C(\mathcal{B})$ spanned by the functions
\[
\psi_n(r)S^l_m(\varphi,\vartheta_1,\dots,\vartheta_{N-2}),
\]
for $1\leq n\leq n_0$, $0\leq m\leq m_0$, and $1\leq l\leq
h(m,N)$. All norms on $E$ are equivalent. In particular, there
exists a constant $C_1=C_1(m_0,n_0)$ such that the
$\varepsilon/3$-neighbourhood of Strassen's ball $\mathcal{K}$ in
the space $E$ equipped by the uniform norm contains in the ball of
radius $1+C_1\varepsilon$ with respect to Strassen's norm. Then we
have
\begin{align*}
\mathsf{P}\{A_{1q}(\varepsilon)\}&\leq\mathsf{P}\left\{\left\|
\frac{\eta^{(m_0,n_0)}_{\mathbf{y}_q,u_q}}{\sqrt{h(t_q)}}\right\|^2_S>(1+
C_1\varepsilon)^2\right\}\\
&=\mathsf{P}\left\{\left\|\eta^{(m_0,n_0)}_{\mathbf{y}_q,u_q}\right\|^2_S>(1+
C_1\varepsilon)^2h(t_q)\right\}\\
&=\mathsf{P}\left\{\sum^{m_0}_{m=0}\sum^{h(m,N)}_{l=1}\sum^{n_0}_{n=1}
\lambda_{mn}(\xi^{lq}_{mn})^2>2(1+C_1\varepsilon)^2h(t_q)\right\}.
\end{align*}
Denote
\[
\lambda=\min_{\substack{0\leq m\leq m_0\\ 1\leq n\leq
n_0}}\lambda_{mn},\qquad M=n_0\sum^{m_0}_{m=0}h(m,N)
\]
and let $\chi_M$ denotes a random variable that has $\chi^2$
distribution with $M$ degrees of freedom. Using standard
probability estimates for $\chi_M$, we can write
\begin{align*}
\mathsf{P}\{A_{1q}(\varepsilon)\}&\leq\mathsf{P}\{\chi_M>2(1+C_1\varepsilon)^2
\lambda^{-1}h(t_q)\}\\
&\leq\exp\{-(1+C_1\varepsilon)^2\lambda^{-1}h(t_q)\}
\end{align*}
for large enough $h(t_q)$. Applying Lemma~\ref{l6},
condition~\ref{l6item2} concludes proof.

Now we prove (\ref{inequal}a). Denote
\begin{align*}
\zeta_1&=\sup_{\substack{(\mathbf{y}_1,u_1)\in R_r(\mathbf{y},u)\\
(\mathbf{y}_2,u_2)\in
R_r(\mathbf{y},u)}}\frac{|\xi(\mathbf{y}_1)-\xi(\mathbf{y}_2)|}
{\sqrt{2}u^H_1},\\
\zeta_2&=\sup_{\substack{(\mathbf{y}_1,u_1)\in R_r(\mathbf{y},u)\\
(\mathbf{y}_2,u_2)\in
R_r(\mathbf{y},u)}}\frac{\displaystyle\sup_{\mathbf{x}\in\mathcal{B}}
|\xi(\mathbf{y}_1+u_1\mathbf{x})-\xi(\mathbf{y}_2+u_2\mathbf{x})|}
{\sqrt{2}u^H_1},\\
\zeta_3&=\sup_{\substack{(\mathbf{y}_1,u_1)\in R_r(\mathbf{y},u)\\
(\mathbf{y}_2,u_2)\in
R_r(\mathbf{y},u)}}\left|\frac{1}{\sqrt{2}u^H_1}-\frac{1}{\sqrt{2}u^H_2}\right|
\sup_{\mathbf{x}\in\mathcal{B}}
|\xi(\mathbf{y}_2+u_2\mathbf{x})-\xi(\mathbf{y}_2)|.
\end{align*}
It is easy to see that $\zeta_1\leq\zeta_2$ and
$\|\eta_{\mathbf{y}_q,u_q}-
\eta_{\mathbf{y},u}\|_{\infty}\leq\zeta_1+\zeta_2+\zeta_3$. It
follows that
\begin{equation}\label{zetaest}
\mathsf{P}\{A_{3q}(\varepsilon)\}\leq
\mathsf{P}\left\{\zeta_2>\varepsilon\sqrt{h(t_q)}/12\right\}
+\mathsf{P}\left\{\zeta_3>\varepsilon\sqrt{h(t_q)}/6\right\}.
\end{equation}
The second term in the right hand side may be estimated as
\begin{equation}\label{zeta3est}
\begin{aligned}
\mathsf{P}\left\{\zeta_3>\varepsilon\sqrt{h(t_q)}/6\right\}&=
\mathsf{P}\left\{\sup_{\substack{\mathbf{x}\in\mathcal{B}\\
(\mathbf{y}_1,u_1)\in R_r(\mathbf{y},u)\\
(\mathbf{y}_2,u_2)\in R_r(\mathbf{y},u)}}|\xi(x)|[(u_2/u_1)^H-1]
>\frac{\varepsilon\sqrt{2h(t_q)}}{6}\right\}\\
&\leq\mathsf{P}\left\{\sup_{\mathbf{x}\in\mathcal{B}}\xi(\mathbf{x})>
\frac{\varepsilon\sqrt{2h(t_q)}}{6\delta}\right\},
\end{aligned}
\end{equation}
where $\delta=\delta(r)=\max\{r,e^{2Hr}-1,e^r-1\}$.

Another large deviation estimate \cite{Lif}*{Section~14, (12)}
states that there exists a constant $C=C(H)$ such that for all
$K>0$
\begin{equation}\label{lif}
\mathsf{P}\left\{\sup_{\mathbf{x}\in\mathcal{B}}\xi(x)>K\right\}\leq
CK^{N/H-1}\exp(-K^2/2).
\end{equation}
Using this fact, we can continue estimate \eqref{zeta3est} as
follows
\[
\mathsf{P}\left\{\zeta_3>\varepsilon\sqrt{h(t_q)}/6\right\}\leq
C\frac{\varepsilon^{N/H-1}[h(t_q)]^{(N-H)/(2H)}}{\delta^{N/H-1}}
\exp\left(-\frac{2\varepsilon^2h(t_q)}{72\delta^2}\right).
\]
If we choose such a small $r$ that $\delta<\varepsilon/6$, then by
Lemma~\ref{l6}, condition~\ref{l6item2} the series
\[
\sum^{\infty}_{q=1}\mathsf{P}\left\{\zeta_3>\varepsilon\sqrt{h(t_q)}/6\right\}
\]
converges.

Using \eqref{prop2}, we write random variable $\zeta_2$ as follows
\[
\zeta_2=\sup_{\substack{\mathbf{x}\in\mathcal{B}\\
(\mathbf{y}_1,u_1)\in R_r(\mathbf{y},u)\\
(\mathbf{y}_2,u_2)\in
R_r(\mathbf{y},u)}}\left|\xi\left(\frac{\mathbf{y}_1+u_1\mathbf{x}}{2^{1/(2H)}u_1}\right)-
\xi\left(\frac{\mathbf{y}_2+u_2\mathbf{x}}{2^{1/(2H)}u_1}\right)\right|.
\]
The right hand side can be estimated as follows.
\begin{align*}
\frac{\|(\mathbf{y}_1+u_1\mathbf{x})-(\mathbf{y}_2+u_2\mathbf{x})\|}{2^{1/(2H)}u_1}&\leq
\frac{\|\mathbf{y}_1-\mathbf{y}_2\|}{2^{1/(2H)}u_1}+\frac{\|\mathbf{x}\|\cdot
|u_1-u_2|}{2^{1/(2H)}u_1}\\
&\leq\frac{2ru}{2^{1/(2H)}u_1}+2^{-1/(2H)}\left(1-\frac{u_2}{u_1}\right)\\
&\leq 2^{1-1/(2H)}re^r+2^{-1/(2H)}(e^{2r}-1)\\
&\leq 2^{1-1/(2H)}\delta e^r+2^{-1/(2H)}(\delta^2+2\delta)\\
&\leq 2^{1-1/(2H)}(\delta^2+\delta)+2^{-1/H}(\delta^2+\delta)\\
&\leq C_2\delta.
\end{align*}

Let $\mathbf{z}_j$, $1\leq j\leq C(C_2\delta)^{-N}$ be the
$C_1\delta$-net in $\mathcal{B}$. Standard entropy estimate for
the first term in the right hand side of \eqref{zetaest} gives
\begin{align*}
\mathsf{P}\left\{\zeta_2>\varepsilon\sqrt{h(t_q)}/12\right\}&\leq
\mathsf{P}\left\{\sup_{\substack{\mathbf{x}\in\mathcal{B}\\
\|\mathbf{y}\|\leq
C_2\delta}}|\xi(\mathbf{x}+\mathbf{y})-\xi(\mathbf{x})|>
\frac{\varepsilon\sqrt{h(t_q)}}{12}\right\}\\
&\leq\mathsf{P}\left\{\sup_{\substack{1\leq j\leq C(C_2\delta)^{-N}\\
\|\mathbf{y}\|\leq
C_2\delta}}|\xi(\mathbf{z}_j+\mathbf{y})-\xi(\mathbf{z}_j)|>
\frac{\varepsilon\sqrt{h(t_q)}}{24}\right\}\\
&\leq C\delta^{-N}\mathsf{P}\left\{\sup_{\|\mathbf{y}\|\leq
1}|\xi(\mathbf{y}|>\frac{\varepsilon\sqrt{h(t_q)}}{24(C_2\delta)^H}\right\}\\
&\leq
C\delta^{-N}\cdot\frac{\varepsilon^{N/H-1}[h(t_q)]^{(N-H)/(2H)}}
{\delta^{N-H}}\exp\left(-\frac{\varepsilon^2h(t_q)}{1152(C_2\delta)^{2H}}\right).
\end{align*}
Here we used Lemma~\ref{l1} and \eqref{lif}. If we choose such a
small $r$ that
\[
\frac{\varepsilon^2}{1152(C_2\delta)^{2H}}>1,
\]
then by Lemma~\ref{l6}, condition~\ref{l6item2} the series
\[
\sum^{\infty}_{q=1}\mathsf{P}\left\{\zeta_2>\varepsilon\sqrt{h(t_q)}/12\right\}
\]
converges. This concludes proof of Lemma~\ref{l3}.

\subsection{Proof of Lemma~\ref{l4}}

By Lemma~\ref{l6}, condition~\ref{l6item1}, for any $q\geq 0$
there exists a number $t'_q\in[t_q,t_{q+1})$ such that
$\mathcal{P}(t)\cap R_{\mathbf{k}_qp_q}\neq\varnothing$. Choose
arbitrary points $(\mathbf{y}'_q,u'_q)\in \mathcal{P}(t)\cap
R_{\mathbf{k}_qp_q}$. It is easy to see that the sequence $t'_q$
satisfies Lemma~\ref{l6}, condition~\ref{l6item2} as well.

The set of all $f\in C(\mathcal{B})$ with $0<\|f\|_S<1$ is dense
in $\mathcal{K}$. It is enough to prove that for any such $f$ we
have
\[
\liminf_{q\to\infty}\left\|\frac{\xi(\mathbf{y}'_q+u'_q\mathbf{x})-\xi(\mathbf{y}'_q)}
{(u'_q)^H\sqrt{2h(t'_q)}}-f\right\|_{\infty}=0\quad\mathsf{P}-\text{a.s.}
\]

According to \cite{Li}*{Theorem~5.1}, there exists a constant
$C_3=C_3(N,H)$ such that for all $\varepsilon\in(0,1]$
\[
\mathsf{P}\left\{\sup_{\mathbf{x}\in\mathcal{B}}|\xi(\mathbf{x})|\leq
\varepsilon\right\}\geq\exp\left(-C_3\varepsilon^{-N/H}\right).
\]
Denote $\beta=(C_3)^{H/N}(1-\|f\|^2_S)^{-H/N}$. We will prove that
\[
\liminf_{q\to\infty}[h(t'_q)]^{(N+2H)/(2N)}\left\|\frac{\xi(\mathbf{y}'_q+u'_q\mathbf{x})-
\xi(\mathbf{y}'_q)}{(u'_q)^H\sqrt{2h(t'_q)}}-f\right\|_{\infty}\leq\frac{1}{\sqrt{2}}\beta
\quad\mathsf{P}-\text{a.s.}
\]

Consider the event
\[
\tilde{A}_{1q}(\varepsilon)=\left\{\left\|\frac{\xi(\mathbf{y}'_q+u'_q\mathbf{x})-
\xi(\mathbf{y}'_q)}{(u'_q)^H}-\sqrt{2h(t'_q)}f\right\|_{\infty}\leq\beta
(1+\varepsilon)[h(t'_q)]^{-H/N}\right\}.
\]
Using \cite{Mon}*{Proposition~4.2}, we obtain
\begin{align*}
\log\mathsf{P}\{\tilde{A}_{1q}(\varepsilon)\}&\geq
2h(t'_q)(-1/2)\|f\|^2_S-C_3\beta^{-N/H}(1+\varepsilon)^{-N/H}h(t'_q)\\
&=-h(t'_q)\left[\|f\|^2_S+(1-\|f\|^2_S)(1+\varepsilon)^{-N/H}\right].
\end{align*}
The multiplier in square brackets is less than $1$. It follows
that
\begin{equation}\label{tildea1}
\sum^{\infty}_{q=0}\mathsf{P}\{\tilde{A}_{1q}(\varepsilon)\}=\infty.
\end{equation}
If the events $\tilde{A}_{1q}(\varepsilon)$ were independent,
the usage of the second Borel--Cantelli lemma would conclude the
proof. However, they are dependent.

In order to create independence, we use another spectral
representation of the multiparameter fractional Brownian motion,
as \cites{Mon,Li} did. Let $\overline{W}$ denotes a complex-valued
scattered Gaussian random measure on $\mathbb{R}^N$ with Lebesgue
measure as its control measure.

\begin{lemma}[Global spectral representation]\label{l7}
There exists a constant $C_4=C_4(N,H)$ such that
\[
\xi(\mathbf{x})=C_4\int_{\mathbb{R}^N}\left(e^{i(\mathbf{p},\mathbf{x})}-1\right)
\|\mathbf{p}\|^{-(N/2)-H}\,d\overline{W}(\mathbf{p}).
\]
\end{lemma}

This result is well-known. Using formulas 2.2.3.1, 2.5.6.1, and
2.5.3.13 from \cite{Pru}, one can prove that
\[
C_4=2^H\sqrt{\frac{H\Gamma((N+H)/2)}{\Gamma(N/2)\Gamma(1-H)}}.
\]

Let $0<a<b$ be two real numbers. Denote
\[
\xi^{(a,b)}(\mathbf{x})=C_4\int_{\|\mathbf{p}\|\in(a,b]}\left(
e^{i(\mathbf{p},\mathbf{x})}-1\right)
\|\mathbf{p}\|^{-(N/2)-H}\,d\overline{W}(\mathbf{p}),\qquad
\tilde{\xi}^{(a,b)}(\mathbf{x})=\xi(\mathbf{x})-\xi^{(a,b)}(\mathbf{x}).
\]

\begin{lemma}\label{l8}
The random field $\tilde{\xi}^{(a,b)}(\mathbf{x})$ has the next
properties.
\begin{enumerate}
\item\label{pro81} It has homogeneous increments.
\item\label{pro82} For any $u\in\mathbb{R}$
\begin{equation}\label{prop82}
\tilde{\xi}^{(a,b)}(u\mathbf{x})\stackrel{d}{=}u^H\tilde{\xi}^{(ua,ub)}(\mathbf{x}).
\end{equation}
\item\label{pro83} It is isotropic.
\end{enumerate}
\end{lemma}

This Lemma can be proved exactly in the same way, as
Lemma~\ref{l1}.

Put
\[
d_q=(u'_q)^{-1}\exp\{h(t'_q)[\exp(h(t'_q)+1-H]\}
\]
and consider the events
\begin{align*}
\tilde{A}_{2q}(\varepsilon)&=\left\{\left\|\frac{\xi^{(d_{q-1},d_q)}(\mathbf{y}'_q+
u'_q\mathbf{x})-\xi^{(d_{q-1},d_q)}(\mathbf{y}'_q)}{(u'_q)^H}-
\sqrt{2h(t'_q)}f\right\|_{\infty}\leq\beta(1+\varepsilon)[h(t'_q)]^{-H/N}\right\},\\
\tilde{A}_{3q}(\varepsilon)&=\left\{\left\|\frac{\tilde{\xi}^{(d_{q-1},d_q)}(\mathbf{y}'_q+
u'_q\mathbf{x})-\tilde{\xi}^{(d_{q-1},d_q)}(\mathbf{y}'_q)}{(u'_q)^H}
\right\|_{\infty}\geq\varepsilon\beta[h(t'_q)]^{-H/N}\right\}.
\end{align*}

\begin{lemma}\label{l9}
We have
\[
\sum^{\infty}_{q=0}\mathsf{P}\{\tilde{A}_{3q}(\varepsilon)\}<\infty.
\]
\end{lemma}

\begin{proof}
Using Lemma~\ref{l8}, one can write
\[
\mathsf{P}\{\tilde{A}_{3q}(\varepsilon)\}=\mathsf{P}\left\{
\|\tilde{\xi}^{(u'_qd_{q-1},u'_qd_q)}(\mathbf{x})\|_{\infty}\geq\varepsilon\beta
[h(t'_q)]^{-H/N}\right\}.
\]
Put
\[
x_q=\exp\{-\exp[(1-\|f\|^2_S)h(t'_q)]\}.
\]
Using Lemma~\ref{l8} once more, we have
\[
\mathsf{P}\{\tilde{A}_{3q}(\varepsilon)\}=\mathsf{P}\left\{
\sup_{\|\mathbf{x}\|\leq
x_q}|\tilde{\xi}^{(x^{-1}_qu'_qd_{q-1},x^{-1}_qu'_qd_q)}
(\mathbf{x})|\geq\varepsilon x^H_q\beta[h(t'_q)]^{-H/N}\right\}.
\]

Denote
\[
\zeta_q(\mathbf{x})=\tilde{\xi}^{(x^{-1}_qu'_qd_{q-1},x^{-1}_qu'_qd_q)}(\mathbf{x}).
\]
We estimate the variance of the random field $\zeta_q(\mathbf{x})$
for $\|\mathbf{x}\|\leq x_q$. We have
\begin{align*}
\mathsf{E}[\zeta_q(\mathbf{x})]^2&=2C^2_4\int_{\|\mathbf{p}\|\leq
x^{-1}_qu'_qd_{q-1}}(1-\cos(\mathbf{p},\mathbf{x}))\|\mathbf{p}\|^{-N-2H}\,d\mathbf{p}\\
&\quad+2C^2_4\int_{\|\mathbf{p}\|>
x^{-1}_qu'_qd_q}(1-\cos(\mathbf{p},\mathbf{x}))\|\mathbf{p}\|^{-N-2H}\,d\mathbf{p}.
\end{align*}
In the first integral, we bound $1-\cos(\mathbf{p},\mathbf{x})$ by
$\|\mathbf{p}\|^2\cdot\|\mathbf{x}\|^2/2$. In the second integral,
we bound it by $2$. Then we have
\[
\mathsf{E}[\zeta_q(\mathbf{x})]^2\leq
C^2_4x^2_q\int_{\|\mathbf{p}\|\leq
x^{-1}_qu'_qd_{q-1}}\|\mathbf{p}\|^{2-N-2H}\,d\mathbf{p}+
4C^2_4\int_{\|\mathbf{p}\|>x^{-1}_qu'_qd_q}\|\mathbf{p}\|^{-N-2H}\,d\mathbf{p}.
\]
Now we pass to spherical coordinates and obtain
\begin{align*}
\mathsf{E}[\tilde{\xi}^{(x^{-1}_qu'_qd_{q-1},x^{-1}_qu'_qd_q)}
(\mathbf{x})]^2&\leq
Cx^2_q\int^{x^{-1}_qu'_qd_{q-1}}_0p^{1-2H}\,dp
+C\int^{\infty}_{x^{-1}_qu'_qd_q}p^{-1-2H}\,dp\\
&=Cx^{2H}_q[(u'_qd_{q-1})^{2-2H}+(u'_qd_q)^{-2H}].
\end{align*}
Substituting definitions of $d_q$ and $x_q$ to the last
inequality, we obtain
\[
\mathsf{E}[\zeta_q(\mathbf{x})]^2\leq
C\exp\left\{-2H\{\exp[(1-\|f\|^2_S)h(t'_q)]+(1-H)h(t'_q)\}\right\}
\]
or,
\[
\mathsf{E}[\tilde{\xi}^{(x^{-1}_qu'_qd_{q-1},x^{-1}_qu'_qd_q)}
(\mathbf{x})-\tilde{\xi}^{(x^{-1}_qu'_qd_{q-1},x^{-1}_qu'_qd_q)}
(\mathbf{y})]^2\leq\varphi^2_q(\|\mathbf{x}-\mathbf{y}\|),
\]
for $\|\mathbf{x}-\mathbf{y}\|=\delta\leq x_q$, where
\[
\varphi^2_q(\delta)=C\min\left\{\delta^{2H},\exp\left\{-2H\{
\exp[(1-\|f\|^2_S)h(t'_q)]+(1-H)h(t'_q)\}\right\}\right\}.
\]

We need the next lemma \cite{Fer}

\begin{lemma}\label{l10}
Let $\zeta(\mathbf{x})$, $\mathbf{x}\in[0,x]^N$ be a separable
centred Gaussian random field. Assume that
\[
\sup_{\substack{\mathbf{x},\mathbf{y}\in[0,x]^N\\\|\mathbf{x}-\mathbf{y}\|\leq\delta}}
\mathsf{E}(\zeta(\mathbf{x})-\zeta(\mathbf{y}))^2\leq\varphi^2(\delta).
\]
Then, for any sequence of positive real numbers $y_0$, $y_1$,
\dots, $y_p$, \dots, and for any sequence of integer numbers
$m_1$, $m_2$, \dots, $m_p$, \dots, every of which can be divided
by previous,
\[
\mathsf{P}\left\{\sup_{\mathbf{x}\in[0,x]^N}|\zeta(\mathbf{x})|\geq
y_0\varphi(x)+\sum^{\infty}_{p=1}y_p\varphi(x/2m_p)\right\}\leq
\sqrt{2/\pi}\sum^{\infty}_{p=0}(m_{p+1})^N\int^{\infty}_{y_p}e^{-u^2/2}\,du.
\]
\end{lemma}

Put $m_{p,q}=q^{2^p}$, $y_{0q}=2\sqrt{(N+1)h(t'_q)}$, and
\[
y_{p,q}=\varepsilon(p+1)^{-2}x^H_q\beta[h(t'_q)]^{-H/N}/\varphi_q(2x_q\cdot
q^{-2^p}),\qquad q\geq 1.
\]
For large enough $q$,
\[
y_{p,q}>2\sqrt{(N+3)h(t'_q)}2^{p/2}
\]
for all $p\geq 1$. Moreover,
\[
y_{0,q}\varphi_q(x_q)+\sum^{\infty}_{p=1}y_{p,q}\varphi_q(x_q/2m_{p,q})<
\varepsilon x^H_q\beta[h(t'_q)]^{-H/N}.
\]
We have
\[
\sum^{\infty}_{q=1}q^2e^{-y^2_{0q}/2}+\sum^{\infty}_{q=1}\sum^{\infty}_{p=1}
q^{N\cdot 2^{p+1}}e^{-y^2_{pq}/2}<\infty,
\]
and application of Lemma~\ref{l10} finishes the proof.
\end{proof}

It follows from definition of the events
$\tilde{A}_{1q}(\varepsilon)$, $\tilde{A}_{2q}(\varepsilon)$, and
$\tilde{A}_{2q}(\varepsilon)$, that
\begin{equation}\label{arelations}
\tilde{A}_{1q}(\varepsilon)\subset\tilde{A}_{2q}(2\varepsilon)\cup
\tilde{A}_{3q}(\varepsilon)\subset\tilde{A}_{1q}(2\varepsilon)\cup
\tilde{A}_{2q}(\varepsilon).
\end{equation}
Combining \eqref{tildea1}, \eqref{arelations}, and Lemma~\ref{l9},
we get
\begin{equation}\label{tildea2}
\sum^{\infty}_{q=0}\mathsf{P}\{\tilde{A}_{2q}(\varepsilon)\}=\infty.
\end{equation}

Now we prove that the events $\tilde{A}_{2q}(\varepsilon)$ are
independent. It is enough to prove that $d_{q-1}<d_q$. Using the
definition of $d_q$, this inequality becomes
\[
u'_q<e^{1-H}u'_{q-1}.
\]
By our choice of $u'_q$, we have $e^{p_qr}\leq u'_q\leq
e^{p_q+1}r$. Since by construction of the parallelepipeds
$R_{\mathbf{k}_qp_q}$ two adjacent parallelepipeds can lie in the
same $u$-layer or in adjacent $u$-layers, we have
\[
u'_{q-1}\leq e^{2r}u'_q.
\]
We choose $r<(1-H)/2$, and we are done.

It follows from the second Borel--Cantelli lemma that
\begin{equation}\label{bc2}
\mathsf{P}\left\{\limsup_{q\to\infty}\tilde{A}_{2q}(\varepsilon)\right\}=1.
\end{equation}
Combining \eqref{arelations}, \eqref{bc2}, and Lemma~\ref{l9}, we
get
\[
\mathsf{P}\left\{\limsup_{q\to\infty}\tilde{A}_{1q}(3\varepsilon)\right\}=1.
\]
Since $\varepsilon$ can be chosen arbitrarily close to $0$,
Lemma~\ref{l4} is proved.

\section{Examples}\label{examples}

\subsection{Local functional law of the iterated logarithm}

Let $t_0=3$. Let $\mathcal{J}(t)$ contains only one element $0$.
Let $\mathbf{y}_0=\mathbf{0}$ and $u_0=t^{-1}$. Then we have
\[
R_1(\mathbf{0},u)=\{\,(\mathbf{y},v)\colon\|\mathbf{y}\|\leq
u,e^{-1}u\leq v\leq eu\,\}.
\]
It is easy to see that $dA_1(u)$ is comparable to
\[
du\int_{\|\mathbf{y}\|\leq
u}\frac{d\mathbf{y}}{u^{N+1}}=C\frac{du}{u}.
\]
The function $h(u)=\log\log u$ satisfies the conditions of
Theorem~\ref{th1}. We obtain, that, in the uniform topology, the
set of $\mathsf{P}$-a. s. limit points of the cloud of increments
\[
\frac{\xi(t\mathbf{x})}{\sqrt{2\log\log t^{-1}}t^H}
\]
as $t\downarrow 0$ is Strassen's ball~$\mathcal{K}$. For the case
of $N=1$ and $H=1/2$, this result is due to \cite{Gan}.

Let $F(f)=\|f\|_{\infty}$, $f\in C(\mathcal{B})$. On the one hand,
we have
\[
\limsup_{t\downarrow 0}\frac{\xi(t)}{\sqrt{2\log\log
t^{-1}}t^H}=\sup_{f\in\mathcal{K}}\|f\|\qquad\mathsf{P}-\text{a.s.}
\]
On the other hand, according to \cite{BeJaRo}, we have
\[
\limsup_{t\downarrow 0}\frac{\xi(t)}{\sqrt{2\log\log
t^{-1}}t^H}=1\qquad\mathsf{P}-\text{a.s.}
\]
It follows that
\begin{equation}\label{supf}
\sup_{f\in\mathcal{K}}\|f\|=1.
\end{equation}

\subsection{Global functional law of the iterated logarithm}

Let $t_0=3$. Let $\mathcal{J}(t)$ contains only one element $0$.
Let $\mathbf{y}_0=\mathbf{0}$ and $u_0=t$. It is easy to check,
that $dA_1(t)$ is comparable to $t^{-1}\,dt$. It follows that, in
the uniform topology, the set of $\mathsf{P}$-a. s. limit points
of the cloud of increments
\[
\frac{\xi(t\mathbf{x})}{\sqrt{2\log\log t}t^H}
\]
as $t\to\infty$ is Strassen's ball~$\mathcal{K}$. For the case of
$N=1$ and $H=1/2$, this result is due to \cite{Str}. Using the
continuous functional $F(f)=\|f\|_{\infty}$, we obtain
\[
\lim_{t\to\infty}\sup_{\|\mathbf{x}\|\leq
t}\frac{\xi(\mathbf{x})}{\sqrt{2\log\log
t}t^H}=\sup_{f\in\mathcal{K}}\|f\|\qquad\mathsf{P}-\text{a.s.}
\]
or, by \eqref{supf},
\[
\lim_{t\to\infty}\sup_{\|\mathbf{x}\|\leq
t}\frac{\xi(\mathbf{x})}{\sqrt{2\log\log
t}t^H}=1\qquad\mathsf{P}-\text{a.s.}
\]

\subsection{Functional L\'{e}vy modulus of continuity}

Let $t_0=2$. Let
$\mathcal{J}(t)=\{\,\mathbf{y}\in\mathbb{R}^N\colon\|\mathbf{y}\|\leq
1-t^{-1}\,\}$ and $u=t^{-1}$ for any
$\mathbf{y}\in\mathcal{J}(t)$. Then we have
\[
\mathcal{P}(t)=\{\,(\mathbf{y},u)\colon\|\mathbf{y}\|\leq
1-t^{-1},u=t^{-1}\,\}
\]
and
\[
\cup_{t\leq
u}\mathcal{P}(t)=\{\,(\mathbf{y},v)\colon\|\mathbf{y}\|\leq
1-u^{-1},u^{-1}\leq v\leq 1\,\}.
\]
It is easy to see that $dA_1(u)$ is comparable to
\[
(1-u^{-1})^N\int^{u^{-1}}_1v^{-N-1}\,dv\sim u^N.
\]
The function $h(u)=N\log u$ satisfies the conditions of
Theorem~\ref{th1}. It follows that, in the uniform topology, the
set of $\mathsf{P}$-a. s. limit points of the cloud of increments
\[
\mathcal{S}(t)=\left\{\,\eta(\mathbf{x})=\frac{\xi(\mathbf{y}+t\mathbf{x})
-\xi(\mathbf{y})}{\sqrt{2N\log t^{-1}}t^H}\colon\|\mathbf{y}\|\leq
1-t\,\right\}
\]
as $t\downarrow 0$ is Strassen's ball~$\mathcal{K}$. For the case
of $N=1$ and $H=1/2$, this result is due to \cite{Mue}. Using the
continuous functional $F(f)=\|f\|_{\infty}$ and \eqref{supf}, we
obtain:
\begin{equation}\label{lmcN}
\limsup_{\|\mathbf{y}\|\downarrow
0}\sup_{\mathbf{x}\in\mathcal{B}}\frac{|\xi(\mathbf{x}+\mathbf{y})-\xi(\mathbf{x})|}
{\sqrt{2N\log\|\mathbf{y}\|^{-1}}\|\mathbf{y}\|^H}=1\qquad\mathsf{P}-\text{a.s.},
\end{equation}
which coincides with the results by \cite{BeJaRo}.

Let $L$ be the linear space of all deterministic functions $f\in
C(\mathcal{B})$ satisfying the condition
\begin{equation}\label{lmcd}
\limsup_{\|\mathbf{y}\|\downarrow
0}\sup_{\mathbf{x}\in\mathcal{B}}\frac{|f(\mathbf{x}+\mathbf{y})-f(\mathbf{x})|}
{\sqrt{2N\log\|\mathbf{y}\|^{-1}}\|\mathbf{y}\|^H}=1.
\end{equation}
It follows from \eqref{lmcN} that $\mu(L)=1$. By
\cite{Lif}*{Section~9, Proposition~1} $\mathcal{H}_{\xi}\subset
L$. From Lemma~\ref{l2} we obtain the following Bernstein-type
theorem from approximation theory:
\begin{theorem}
Let $f\in C(\mathcal{B})$ with $f(\mathbf{0})=0$ satisfies the
condition
\[
\sum^{\infty}_{m=0}\sum^{\infty}_{n=1}\sum^{h(m,N)}_{l=1}
\frac{(f^l_{mn})^2}{\lambda_{mn}}<\infty.
\]
Then $f$ satisfies \eqref{lmcd}.
\end{theorem}

\end{document}